\documentstyle{article}

\newcommand{\N}{{\mbox{\bf N}}}

\def\cirk{\,{\raisebox{.3ex}{\tiny $\circ$}}\,}

\def\prop#1#2{\vspace{2ex} \noindent{\sc #1.} {\it #2} \par \vspace{2ex}}

\def\Mat{{\mbox{${\mbox{\bf Mat}}_{\cal F}$}}}

\def\Vect{{\mbox{${\mbox{\bf Vect}}_{\cal F}$}}}

\def\Vectfd{{\mbox{${\mbox{\bf Vect}}_{\cal F}^{\mbox{\scriptsize \it fd}}$}}}

\begin{document}

\title
{{\sc Symmetric Self-Adjunctions:}\\A Justification of Brauer's
Representation of Brauer's Algebras}

\author{
{\sc Kosta Do\v sen} and {\sc Zoran Petri\' c}
\\[.3cm]{\small Mathematical Institute, SANU}
\\[-.1cm]{\small Knez Mihailova 35, P.O. Box 367}
\\[-.1cm]{\small 11001 Belgrade, Serbia}
\\[-.1cm]{\small email: \{kosta, zpetric\}@mi.sanu.ac.yu}}

\date{ }
\maketitle

\begin{abstract}
\noindent A classic result of representation theory is Brauer's
construction of a diagrammatical (geometrical) algebra whose
matrix representation is a certain given matrix algebra, which is
the commutating algebra of the enveloping algebra of the
representation of the orthogonal group. The purpose of this paper
is to provide a motivation for this result through the categorial
notion of symmetric self-adjunction.
\end{abstract}

\vspace{.3cm}

{\small

\vspace{.5ex}

\noindent {\it Keywords$\,$}: Brauer's centralizer algebras,
matrix representation, orthogonal group, adjoint functor

\vspace{.5ex}

\noindent {\it Acknowledgement$\,$}: We are grateful to Professor
Neda Bokan for enabling us to present the results of this paper at
the conference \emph{Contemporary Geometry and Related Topics}
(Belgrade, June 2005). The writing of this paper was supported by
the Ministry of Science of Serbia (Grant 144013, Representation of
Logical Structures). \vspace{1cm} }

\section{Introduction}

In \cite{B37} Richard Brauer introduced a class of diagrammatical
algebras and found a matrix representation for them. These
algebras have been in focus again after more than fifty years
since Temperley-Lieb algebras, the subalgebras of Brauer's
algebras, have started to play an important role in knot theory
and low-dimensional topology via the polynomial approach to knot
invariants (see \cite{KL94}, \cite{L97} and \cite{PS97}). Brauer's
algebras arose as a side product of investigations concerning a
representation of the orthogonal group ${\cal O}(n)$, and the
representation of these algebras was established after Brauer's
remark that if one associates (in a natural way) diagrams to
matrices from a particular class, then the product of these
matrices corresponds to the operation of ``composition'' of such
diagrams. We give more details concerning Brauer's introduction of
these algebras in Section~2. So, for a given matrix algebra,
Brauer had constructed a diagrammatical algebra whose matrix
representation turned out to be this matrix algebra. We find such
a representation insufficiently justified.

On the other hand, there are several results concerning
diagrammatical characterization of various kinds of free
adjunctions and related notions from category theory (see for
example \cite{D99}, \cite{DP03a}, \cite{DP03b}, \cite{DP05a} and
\cite{KM71}, \cite{FY92}, \cite{JS93}). Combining the fact that
the symmetric self-adjunction freely generated by a singleton set
of objects is isomorphic to the category of Brauer's diagrams and
the fact that a symmetric self-adjunction exists in the skeleton
of the category of finite dimensional vector spaces over a field
$\cal F$, one can find a matrix representation of Brauer's
diagrams that coincides with Brauer's representation. We find this
is a natural justification of Brauer's representation.

\section{Brauer's algebras and their representation}

For every $n\in \N^+$, Brauer's algebra $B_n$ over a field $\cal
F$ of characteristic 0 is a vector space whose basis consists of
$(2n-1)!!$ diagrams, which we call \emph{Brauer's}
$n$-\emph{diagrams} or just $n$-\emph{diagrams}. Every $n$-diagram
consists of $n$ \emph{vertices} in the \emph{top row} and $n$
vertices in the \emph{bottom row}. Each of these $2n$ vertices is
connected by a \emph{thread} with exactly one of the remaining
$2n-1$ vertices. For example,

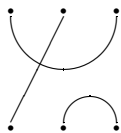
\begin{center}
\begin{picture}(40,50)

\put(0,5){\line(1,2){20}}

\put(0,3){\circle*{2}} \put(0,47){\circle*{2}}
\put(20,3){\circle*{2}} \put(20,47){\circle*{2}}
\put(40,3){\circle*{2}} \put(40,47){\circle*{2}}

\put(30,5){\oval(20,20)[t]} \put(20,45){\oval(40,40)[b]}

\end{picture}
\end{center}

\noindent is a 3-diagram.

So, addition and multiplication by scalars is formal in $B_n$ and
as a vector space, $B_n$ is isomorphic to ${\cal F}^{(2n-1)!!}$.
For the structure of algebra in $B_n$ it is sufficient to define
multiplication of $n$-diagrams. (We call this multiplication
\emph{composition} and denote it by $\cirk$.) To define the
$n$-diagram $D_2\cirk D_1$ for two $n$-diagrams $D_1$ and $D_2$,
we have to identify the bottom row of $D_1$ with the top row of
$D_2$ so that the top row of $D_1$ becomes the top row of
$D_2\cirk D_1$ and the bottom row of $D_2$ becomes the bottom row
of $D_2\cirk D_1$. The threads of $D_2\cirk D_1$ are obtained by
concatenating the threads of $D_1$ and $D_2$. The number $k\geq 0$
of \emph{circular components} that may occur in this procedure
reflects in the scalar $p^k$ which multiplies the resulting
$n$-diagram ($p$ is here a fixed positive integer and the choice
to represent a circle in a diagram by multiplying the rest of the
diagram by $p$ is forced by the matrix algebra in which Brauer
represented $B_n$). For example, let $D_1$ and $D_2$ be the
following 3-diagrams:

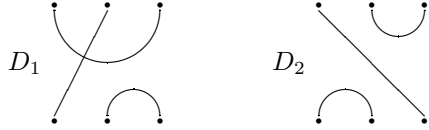
\begin{center}
\begin{picture}(140,50)

\put(0,5){\line(1,2){20}}

\put(0,3){\circle*{2}} \put(0,47){\circle*{2}}
\put(20,3){\circle*{2}} \put(20,47){\circle*{2}}
\put(40,3){\circle*{2}} \put(40,47){\circle*{2}}

\put(30,5){\oval(20,20)[t]} \put(20,45){\oval(40,40)[b]}

\put(-5,25){\makebox(0,0)[r]{$D_1$}}

\put(140,5){\line(-1,1){40}}

\put(100,3){\circle*{2}} \put(100,47){\circle*{2}}
\put(120,3){\circle*{2}} \put(120,47){\circle*{2}}
\put(140,3){\circle*{2}} \put(140,47){\circle*{2}}

\put(110,5){\oval(20,20)[t]} \put(130,45){\oval(20,20)[b]}

\put(95,25){\makebox(0,0)[r]{$D_2$}}

\end{picture}
\end{center}

\noindent After identification of the bottom row of $D_1$ with the
top row of $D_2$ we have

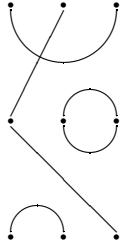
\begin{center}
\begin{picture}(40,90)

\put(0,49){\line(1,2){20}}

\put(0,91){\circle*{2}} \put(20,91){\circle*{2}}
\put(40,91){\circle*{2}}

\put(30,49){\oval(20,20)[t]} \put(20,89){\oval(40,40)[b]}

\put(40,5){\line(-1,1){40}}

\put(0,3){\circle*{2}} \put(0,47){\circle*{2}}
\put(20,3){\circle*{2}} \put(20,47){\circle*{2}}
\put(40,3){\circle*{2}} \put(40,47){\circle*{2}}

\put(10,5){\oval(20,20)[t]} \put(30,45){\oval(20,20)[b]}

\end{picture}
\end{center}

\noindent and $D_2\cirk D_1$ is the following element of $B_3$

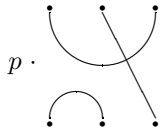
\begin{center}
\begin{picture}(40,50)

\put(40,5){\line(-1,2){20}}

\put(0,3){\circle*{2}} \put(0,47){\circle*{2}}
\put(20,3){\circle*{2}} \put(20,47){\circle*{2}}
\put(40,3){\circle*{2}} \put(40,47){\circle*{2}}

\put(10,5){\oval(20,20)[t]} \put(20,45){\oval(40,40)[b]}

\put(-5,25){\makebox(0,0)[r]{$p\;\cdot$}}

\end{picture}
\end{center}

Let $I_n$ be the $n$-diagram in which for every
${i\in\{1,\ldots,n\}}$ we have that the $i$-th vertex from the top
row is connected with the $i$-th vertex from the bottom row. For
example $I_3$ is

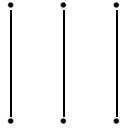
\begin{center}
\begin{picture}(40,50)

\put(0,5){\line(0,1){40}} \put(20,5){\line(0,1){40}}
\put(40,5){\line(0,1){40}}

\put(0,3){\circle*{2}} \put(0,47){\circle*{2}}
\put(20,3){\circle*{2}} \put(20,47){\circle*{2}}
\put(40,3){\circle*{2}} \put(40,47){\circle*{2}}

\end{picture}
\end{center}

\noindent It is pretty obvious that $I_n$ is the unit for
composition and associativity of composition is straightforward
when we rely on such an informal (pictorial) definition of
$\cirk$. For a formal proof of associativity of $\cirk$ one may
consult \cite{DP03c} and \cite{DP03d}. We explain below how these
algebras arose in the work of Brauer.

Let $\cal F$ be a field of characteristic 0 and let $\cal G$ be
${\cal O}(p)$ (group of orthogonal linear transformations of the
$p$-dimensional vector space ${\cal F}^p$ over $\cal F$). Every
member of $\cal G$ is given by an orthogonal $p\times p$ matrix
$G$ ($G^{-1}=G^T$) with entries from $\cal F$.

Brauer was particulary interested in the following representation
of $\cal G$:

\[M(G)=G^{\otimes n}=\underbrace{G\otimes\ldots\otimes G}_{n
\mbox{\scriptsize\rm-fold Kronecker product}} \in \mbox{\rm End}\,
({\cal F}^{p^n})\]

\noindent (That $M$ is a representation, i.e.\ that
$M(G_1G_2)=M(G_1)M(G_2)$, follows from the functoriality of
$\otimes$.)

Let $\cal M$ be the group ${\{M(G)\mid G\in{\cal G}\}}$ and let
$\cal A$ be the \emph{enveloping algebra} of $\cal M$, i.e.,

\[
{\cal A}=\{c_1M(G_1)+...+c_kM(G_k)\mid k\in\N^+,c_i\in{\cal
F},G_i\in{\cal G},1\leq i\leq k\}.
\]

\noindent Brauer's goal was to characterize elements of this
algebra. For this purpose he used the \emph{commutating algebra}
$\cal B$ of $\cal A$, i.e.,

\[
{\cal B}=\{B\in\mbox{\rm End}\,({\cal F}^{p^n})\mid(\forall
A\in{\cal A})\; AB=BA\}.
\]

\noindent The algebra $\cal B$ is a matrix algebra whose elements
are $p^n\times p^n$ matrices of the form

\[
[b_{i,j}]_{p^n\times p^n}= \left[
\begin{array}{cccc}
  b_{1,1}   & b_{1,2}   & \ldots & b_{1,p^n} \\
  b_{2,1}   & b_{2,2}   & \ldots & b_{2,p^n} \\
  \ldots    & \ldots    & \ldots & \ldots \\
  b_{p^n,1} & b_{p^n,2} & \ldots & b_{p^n,p^n} \\
\end{array}
\right]
\]

\noindent with entries from $\cal F$.

We explain now some technical notation that we are going to use
below. There are $p^n$ functions from $\{1,\ldots,n\}$ to
$\{1,\ldots,p\}$. Each of these functions can be envisaged as a
sequence of length $n$ of elements of $\{1,\ldots,p\}$. The set of
these sequences may be ordered lexicographically so that
${(1,1,\ldots,1)}$ is the first and ${(p,\ldots,p)}$ is the last
($p^n$-th) in this ordering. We use this ordering to identify the
elements of $\{1,\ldots,p^n\}$ with the functions from
$\{1,\ldots,n\}$ to $\{1,\ldots,p\}$. So, for
$i\in\{1,\ldots,p^n\}$ and $k\in\{1,\ldots,n\}$, by
$i(k)\in\{1,\ldots,p\}$ we mean the image of $k$ by the $i$-th
function from $\{1,\ldots,n\}$ to $\{1,\ldots,p\}$.

Let

\[
u^i= \left[
\begin{array}{c}
  u^i_1 \\[.2ex]
  \vdots \\[.2ex]
  u^i_p \\
\end{array}
\right] \quad \mbox{\rm and}\quad v^j= \left[
\begin{array}{c}
  v^j_1 \\[.2ex]
  \vdots \\[.2ex]
  v^j_p \\
\end{array}
\right]
\]

\noindent for ${i,j\in\{1\ldots,n\}}$ be $2n$ vectors of ${\cal
F}^p$ and let $\cal G$ acts on them, so that for $G\in{\cal G}$ we
have $(u^i)'=Gu^i$ and $(v^j)'=Gv^j$. Then the following holds.

\prop{Theorem \mbox{\rm (Brauer)}}{The function \[
J(u^1,\ldots,u^n,v^1,\ldots,v^n)=\sum_{1\leq i,j\leq
p^n}b_{i,j}u^1_{i(1)}\ldots u^n_{i(n)}v^1_{j(1)}\ldots
v^n_{j(n)}\] is an invariant of $\cal G$ (i.e.,
$J(u^1,\ldots,v^n)=J((u^1)',\ldots,(v^n)')$) iff the matrix
$[b_{i,j}]_{p^n\times p^n}$ belongs to $\cal B$. }

From the main theorem of invariant theory concerning the
orthogonal group case (see \cite{W46}, Chapter II, Section A.9) it
follows that

\[
J(u^1,\ldots,u^n,v^1,\ldots,v^n)=\sum_{1\leq i,j\leq
p^n}b_{i,j}u^1_{i(1)}\ldots u^n_{i(n)}v^1_{j(1)}\ldots
v^n_{j(n)}\]

\noindent is an invariant of $\cal G$ iff
$J(u^1,\ldots,u^n,v^1,\ldots,v^n)$ is a linear combination of
products of scalar products of the form

\[
(w^1w^2)\cdot(w^3w^4)\cdot\ldots\cdot(w^{2n-1}w^{2n})
\]

\noindent where $w^1w^2\ldots w^{2n}$ is a permutation of vectors
${u^1\ldots u^nv^1\ldots v^n}$.

Up to commutativity of multiplication and scalar product there are
${(2n-1)!!}$ different terms of this form. It is natural to
associate with every such term an $n$-diagram. In this diagram the
vertices from the top row represent the vectors $u^1,\ldots,u^n$
and vertices from the bottom row represent the vectors
$v^1,\ldots,v^n$. Then every thread of the diagram shows which
pairs of vectors occur in the scalar products of the term.
However, this correspondence is a bijection (up to commutativity
of multiplication and scalar product). This means that starting
from an arbitrary $n$-diagram one can find a term of the above
form representing an invariant of $\cal G$ from which a matrix
from $\cal B$ can be extracted. In this way Brauer obtained a
function that maps $n$-diagrams to matrices from $\cal B$.

Brauer's remark was that the result of composition of two
$n$-diagrams is mapped to the product of matrices corresponding to
these diagrams. This is the core of Brauer's representation of
Brauer's algebra since, by linearity, the above correspondence
between $n$-diagrams and matrices can be extended to a
representation of $B_n$ in a unique way.

We illustrate this representation by an example in which ${p=2}$,
${n=3}$ and $D$ is the 3-diagram:

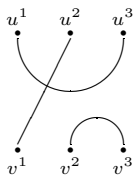
\begin{center}
\begin{picture}(40,50)

\put(0,5){\line(1,2){20}}

\put(0,3){\circle*{2}} \put(0,47){\circle*{2}}
\put(20,3){\circle*{2}} \put(20,47){\circle*{2}}
\put(40,3){\circle*{2}} \put(40,47){\circle*{2}}

\put(30,5){\oval(20,20)[t]} \put(20,45){\oval(40,40)[b]}

\put(0,0){\makebox(0,0)[t]{\scriptsize $v^1$}}
\put(20,0){\makebox(0,0)[t]{\scriptsize $v^2$}}
\put(40,0){\makebox(0,0)[t]{\scriptsize $v^3$}}

\put(0,50){\makebox(0,0)[b]{\scriptsize $u^1$}}
\put(20,50){\makebox(0,0)[b]{\scriptsize $u^2$}}
\put(40,50){\makebox(0,0)[b]{\scriptsize $u^3$}}

\end{picture}
\end{center}

\noindent Then the following term corresponds to $D$

\[
(u^1u^3)\cdot(u^2v^1)\cdot(v^2v^3)=(u^1_1\cdot u^3_1+u^1_2\cdot
u^3_2)\cdot(u^2_1\cdot v^1_1+u^2_2\cdot v^1_2)\cdot(v^2_1\cdot
v^3_1+v^2_2\cdot v^3_2).
\]

\noindent After distributions at the right-hand side of this
equation we obtain a term of the form

\[
\sum_{1\leq i,j\leq 8}b_{i,j}u^1_{i(1)}u^2_{i(2)}u^3_{i(3)}
v^1_{j(1)}v^2_{j(2)}v^3_{j(3)}\]

\noindent where $b_{ij}=1$ iff $i(1)=i(3)$, $i(2)=j(1)$ and
$j(2)=j(3)$, otherwise $b_{ij}=0$. Roughly speaking, $b_{ij}=1$ if
and only if $i$ as a ternary sequence of elements of $\{1,2\}$
above $j$ as a ternary sequence of elements of $\{1,2\}$ as in the
picture below

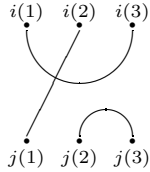
\begin{center}
\begin{picture}(40,50)

\put(0,5){\line(1,2){20}}

\put(0,3){\circle*{2}} \put(0,47){\circle*{2}}
\put(20,3){\circle*{2}} \put(20,47){\circle*{2}}
\put(40,3){\circle*{2}} \put(40,47){\circle*{2}}

\put(30,5){\oval(20,20)[t]} \put(20,45){\oval(40,40)[b]}

\put(0,0){\makebox(0,0)[t]{\scriptsize $j(1)$}}
\put(20,0){\makebox(0,0)[t]{\scriptsize $j(2)$}}
\put(40,0){\makebox(0,0)[t]{\scriptsize $j(3)$}}

\put(0,50){\makebox(0,0)[b]{\scriptsize $i(1)$}}
\put(20,50){\makebox(0,0)[b]{\scriptsize $i(2)$}}
\put(40,50){\makebox(0,0)[b]{\scriptsize $i(3)$}}

\end{picture}
\end{center}

\noindent is ready to ``accept'' $D$, in the sense that linked
elements of $\{1,2\}$ are equal.

The Temperley-Lieb algebra $TL_n$ is a subalgebra of $B_n$ whose
basis consists of non-intersecting $n$ diagrams. The number of
such diagrams is $(2n)!/(n!(n+1)!)$ (the $n$-th Catalan number).
It is proved in \cite{J94} (see also \cite{DKP}) that the
restriction to $TL_n$ of Brauer's representation of $B_n$ is
faithful for $p\geq 2$, which means that this representation is an
embedding of Temperley-Lieb algebras into ${\mbox{\rm End}\,({\cal
F}^{p^n})}$. However, this cannot always be the case for Brauer's
representation of $B_n$ since $(2n-1)!!$ as the dimension of $B_n$
may exceed $p^{2n}$ as the dimension of ${\mbox{\rm End}\,({\cal
F}^{p^n})}$.

We are going now to generalize the notion of Brauer's $n$-diagram
in the sense that we allow different number of vertices in its top
and bottom row. So, let an $m$-$n$-\emph{diagram} be a diagram
like Brauer's $n$-diagram save that it has $m$ vertices in the top
row and $n$ vertices in the bottom row for $m$ not necessarily
equal to $n$. Then we can take instead of just $n$-diagrams for a
particular $n\in\N^+$, the class of $m$-$n$-diagrams for all
$m,n\in\N$ and define the composition of an $m$-$n$-diagram and an
$n$-$q$-diagram analogously to what we had for two $n$-diagrams.
So the result of this composition is an $m$-$q$-diagram multiplied
by a scalar of the form $p^k$ which reflects the number $k$ of
circular components that arise after concatenating the threads of
these diagrams. In this way we obtain the category
$\mbox{\emph{Br}}_p$ whose objects are natural numbers and whose
arrows are $m$-$n$-diagrams with coefficients of the form $p^k$
for fixed $p\geq 1$. In \cite{DP05}, Section 2.3, the category
\emph{Br} related to $\mbox{\emph{Br}}_p$ (case $p=1$) is defined
in a more formal way. One can call this generalization of $B_n$ a
\emph{categorification} of multiplicative submonoids of Brauer's
algebras generated out of the basis. We shall see in the following
section that the category $\mbox{\emph{Br}}_p$ is strongly
connected to the notion of symmetric self-adjunction of
\cite{DP05a} in a way that it may be called the \emph{geometry} of
symmetric self-adjunctions.

All the above shows that Brauer's representation of Brauer's
algebras ``works''. But one may still ask why does it work? Or,
what mathematics underlies this representation? We try to answer
these questions in the following sections.

\section{Symmetric self-adjunctions}

In the hierarchy of categorial notions, one of the topmost
positions is reserved for the notion of adjunction. This notion
can be defined equationally in the following manner: an
\emph{adjunction} is a 6-tuple $\langle {\cal A},{\cal
B},F,G,\varphi,\gamma\rangle$ where

$\cal A$ and $\cal B$ are categories;

$F\!:{\cal B}\rightarrow{\cal A}$ and $G\!:{\cal
A}\rightarrow{\cal B}$ are functors;

$\varphi\!:FG\rightarrow 1_{\cal A}$ and $\gamma\!:1_{\cal
B}\rightarrow GF$ are natural transformations such that the
following \emph{triangular equations} hold in $\cal A$ and $\cal
B$ respectively

\[
\begin{array}{l}
\varphi_{FB}\cirk F\gamma_B=1_{FB},
\\[.2ex]
G\varphi_A\cirk\gamma_{GA}=1_{GA}.
\end{array}
\]

\noindent The definition above is equational in the sense that it
is possible to present the notions of category, functor and
natural transformation equationally. Such an equational definition
guarantees the existence of some free structures on which we rely
in this section.

In \cite{D99}, Section 4.10.1, an $m$-$n$-diagram is associated to
every canonical arrow $f$ of an adjunction. These diagrams are
called the \emph{set of links} of $f$ and are denoted by
$\Lambda(f)$. For example

\begin{center}
\begin{picture}(220,50)

\put(50,20){\oval(12,12)[b]}

\put(35,10){\makebox(0,0)[r]{$\Lambda(\varphi_A)$ is }}

\put(70,0){\makebox(0,0)[r]{$A$}}
\put(70,27){\makebox(0,0)[r]{$F\;G\;A$}}

\put(200,7){\oval(12,12)[t]}

\put(185,10){\makebox(0,0)[r]{and \quad\quad $\Lambda(\gamma_B)$
is }}

\put(220,27){\makebox(0,0)[r]{$B$}}
\put(220,0){\makebox(0,0)[r]{$G\;F\;B$}}

\end{picture}
\end{center}

\noindent and $\Lambda(g\cirk f)=\Lambda(g)\cirk \Lambda(f)$
(where $\cirk$ on the right-hand side denotes the composition of
$m$-$n$-diagrams). We illustrate the soundness of $\Lambda$ with
the first of the triangular equations; this yields the following
picture

\begin{center}
\begin{picture}(200,60)

\put(60,36){\oval(12,12)[t]}

\put(30,15){\makebox(0,0)[r]{$\Lambda(\varphi_{FB})$}}
\put(30,45){\makebox(0,0)[r]{$\Lambda(F\gamma_B)$}}

\put(48,24){\oval(12,12)[b]}

\put(60,60){\makebox(0,0){$F\;B$}}
\put(60,30){\makebox(0,0){$F\;G\;F\;B$}}
\put(60,0){\makebox(0,0){$F\;B$}}

\put(43,36){\line(1,2){9}} \put(54,6){\line(1,2){9}}

\put(160,60){\makebox(0,0){$F\;B$}}
\put(160,0){\makebox(0,0){$F\;B$}}

\put(150,30){\makebox(0,0)[r]{$\Lambda(1_{FB})$}}

\put(154,6){\line(0,1){49}}

\end{picture}
\end{center}

As we have already mentioned, the equational definition of
adjunction guarantees the existence of the adjunction freely
generated by a pair of sets of objects (discrete categories). Then
$\Lambda$ gives rise to functors from both categories involved in
this freely generated adjunction to the category
$\mbox{\emph{Br}}_p$. It is proved in \cite{D99}, Proposition in
Section 4.10.1, that both of these functors are faithful. However,
not all the $m$-$n$-diagrams are covered by the arrows of freely
generated adjunction. It is easy to see that all the
$m$-$n$-diagrams corresponding to these arrows are of
Temperley-Lieb kind (they are non-intersecting diagrams) and even
not all the diagrams of the Temperley-Lieb kind are covered by
this correspondence. This was a motivation for a step leading from
the notion of adjunction to a more specific notion of
self-adjunction (see \cite{DP03a}, \cite{DP03b} and references
therein, see also \cite{DP05a} for a more gradual introduction of
this notion). A \emph{self-adjunction} (also called $\cal
L$-adjunction in \cite{DP03a}) may be introduced as a quadruple
${\langle{\cal A},F,\varphi,\gamma\rangle}$ such that
${\langle{\cal A},{\cal A},F,F,\varphi,\gamma\rangle}$ is an
adjunction. (So, $F\!:{\cal A}\rightarrow{\cal A}$ is an
endofunctor adjoint to itself.)

As in the case of adjunction, one may construct the
self-adjunction freely generated by an arbitrary set of objects.
Then $\Lambda$, as before, gives rise to a functor from the
category involved in this freely generated self-adjunction to the
category $\mbox{\emph{Br}}_p$. This time, the functor is not
faithful because a simple counting of circular components that
occur in compositions of $m$-$n$-diagrams is not sufficient. The
faithfulness of this functor requires some adjustments in the
category $\mbox{\emph{Br}}_p$; namely, one must take into account
not just the number of circular components, but also their
positions in the diagram. (See \cite{DP03a} for the definition of
friezes and $\cal L$-equivalence between them.) It is shown in
\cite{DP03a} that the arrows of freely generated self-adjunction
cover by $\Lambda$ all the diagrams of the Temperley-Lieb kind.

Since all the intersecting $m$-$n$-diagrams are still out of the
range of $\Lambda$ we can make a step forward, to arrive at the
notion of symmetric self-adjunction, which is defined as follows.
A \emph{symmetric self-adjunction} is a quintuple ${\langle{\cal
A},F,\varphi,\gamma,\chi\rangle}$ for which ${\langle{\cal
A},F,\varphi,\gamma\rangle}$ is a self-adjunction, $\chi$ is a
natural transformation from $F\cirk F$ to $F\cirk F$ such that the
equations

\begin{tabbing}
\hspace{2em}\= $\chi_A\cirk\chi_A=1_{FFA}$,\hspace{7em} \=
$\chi_{FA}\cirk F\chi_A\cirk\chi_{FA}=F\chi_A\cirk\chi_{FA}\cirk
F\chi_A$,
\\[.2ex]
\> $\varphi_A\cirk\chi_A=\varphi_A$, \>
$\chi_A\cirk\gamma_A=\gamma_A$,
\\[.2ex]
\> $\varphi_{FA}\cirk F\chi_A=F\varphi_A\cirk \chi_{FA}$, \>
$\chi_{FA}\cirk F\gamma_A=F\chi_A\cirk\gamma_{FA}$.
\end{tabbing}

\noindent are satisfied.

If we extend $\Lambda$ to cover all the canonical arrows of a
symmetric self-adjunction by defining $\Lambda(\chi_A)$ to be

\begin{center}
\begin{picture}(40,30)

\put(20,30){\makebox(0,0){$F\;F\;A$}}
\put(20,0){\makebox(0,0){$F\;F\;A$}}

\put(9,6){\line(1,2){9}} \put(19,6){\line(-1,2){9}}

\end{picture}
\end{center}

\noindent then $\Lambda$ gives rise to a functor from the category
involved in the symmetric self-adjunction freely generated by a
set of objects, to the category $\mbox{\emph{Br}}_p$. It is shown
in \cite{DP05a} that this functor is faithful and, moreover, if
the set of generating objects is a singleton, then this functor is
an isomorphism.

\section{Symmetric self-adjunction of the category \Mat}

In this section we discuss an example of symmetric
self-adjunction. To find such an example we start with the
category \Vect\ of vector spaces over the field $\cal F$. For a
given vector space $V\in\Vect$, the functor
$F\!:\Vect\rightarrow\Vect$ which acts on objects as ${FU=V\otimes
U}$ has the right adjoint ${G\!:\Vect\rightarrow\Vect}$ which maps
a vector space $W$ to the vector space of all linear
transformations from $V$ to $W$.

If we replace \Vect\ by \Vectfd, the compact closed category of
finite dimensional vector spaces over $\cal F$ (for the notion of
compact closed category see \cite{K72}, \cite{KL80} or
\cite{D77}), then the right adjoint of
${F=V\otimes\underline{\;\;\;}}$ becomes
${V^*\otimes\underline{\;\;\;}}$, where $V^*$ is the dual vector
space of $V$. Since  $V$ and $V^*$ are isomorphic in \Vectfd\ we
have that this isomorphism leads to the isomorphism of functors
${V\otimes\underline{\;\;\;}}$ and
${V^*\otimes\underline{\;\;\;}}$. So,
${V\otimes\underline{\;\;\;}}$ becomes a self-adjoint functor. It
is easy to see that the natural transformation $\chi$ (symmetry)
with all the required equations is present in \Vectfd\ with
${F=V\otimes\underline{\;\;\;}}$ and, hence, we obtain an example
of symmetric self-adjunction.

We can simplify the category \Vectfd\ by passing to its skeleton
\Mat\  (a full subcategory of \Vectfd\ such that each object of
\Vectfd\ is isomorphic to exactly one object of \Mat). The
category \Mat\ still provides an example of symmetric
self-adjunction. We can envisage this category as the category
whose objects are natural numbers (the dimensions of finite
dimensional vector spaces) and an arrow ${M\!:m\rightarrow n}$ is
an ${n\times m}$ matrix with entries from $\cal F$. Composition of
such arrows becomes matrix multiplication.

For every $p\in N$, the functor ${p\otimes\underline{\;\;\;}}$,
which maps an object $n$ of \Mat\ to the product $m\cdot n$ and an
arrow $M$ of \Mat\ to the Kronecker product $I_p\otimes M$, is the
part of a symmetric self-adjunction
${\langle\Mat,p\otimes\underline{\;\;\;},\varphi,\gamma,\chi\rangle}$
for some indexed sets $\varphi$, $\gamma$ and $\chi$ of matrices.

Suppose now that $p$ is equal to the $p$ we used in the definition
of composition of $m$-$n$-diagrams. We have the following picture
in which $\cal K$ is the category of the symmetric self-adjunction
freely generated by a single object,

\[
\mbox{$\mbox{\emph{Br}}_p$}\stackrel{\cong}{\longrightarrow}{\cal
K}\stackrel{R}{\longrightarrow}\Mat
\]

\noindent and $R$ is the functor which strictly preserves the
structure of symmetric self-adjunction, and which extends the
function that maps the unique generator of $\cal K$ to the object
1 of \Mat. The functor $R$ exists by the freedom of $\cal K$.

The above composition of functors (which we also denote by $R$)
has the following properties: for every pair of $m$-$n$-diagrams
$D_1$ and $D_2$ we have that ${R(D_2\cirk D_1)=R(D_2)\cdot
R(D_1)}$ (by the functoriality of $R$) and $R(p^k\cdot
D_1)=p^k\cdot R(D_1)$. These properties show that $R$ can serve as
a core for a matrix representation of Brauer's algebras. It is not
difficult to check that this representation coincides with
Brauer's representation, which is now properly justified via
symmetric self-adjunctions.

\end{document}